\documentclass{arXSiH}

\newtheorem{theorem}{Theorem}

\newnumbered{assertion}{Assertion}    
\newnumbered{conjecture}{Conjecture}  
\newnumbered{definition}[theorem]{Definition}
\newnumbered{hypothesis}{Hypothesis}
\newnumbered{remark}[theorem]{Remark}
\newnumbered{construction}[theorem]{Construction}
\newnumbered{constructions}[theorem]{Constructions}
\newnumbered{remarks}[theorem]{Remarks}
\newnumbered{note}{Note}
\newnumbered{notes}[theorem]{Notes}
\newnumbered{observation}{Observation}
\newnumbered{problem}{Problem}
\newnumbered{question}{Question}
\newnumbered{algorithm}{Algorithm}
\newnumbered{example}[theorem]{Example}
\newnumbered{examples}[theorem]{Examples}
\newunnumbered{notation}{Notation} 
\newnumbered{notations}[theorem]{Notations}
\newunnumbered{unexample}{Example}
\newunnumbered{unremark}{Remark}
\newunnumbered{Sconventions}{Some conventions}

\title[A holomorphic map]%
             {A holomorphic map in infinite dimensions}

\author{S. Hiltunen}

\classno{46G20}

\begin{document}
\maketitle

\begin{abstract} We prove, in great detail, holomorphy \œ$E\sqcap
  C\ssp(\spp I\sp,\varPi\sp)\to C\ssp(\spp I\sp,\varPi\sp)$ of the map \œ$
(\sp x\ssp,y\sp)\mapsto x\circ[\,\roman{id\,},y\,]$ where \œ$[\,\roman{id\,},
y\,]:I\owns t\mapsto(\sp t\ssp,y\ssp(t))$ for a real compact interval $I\sp$,
and where $\varPi$ is a complex Banach space and $E$ is a certain locally
convex space of continuous functions \œ$
x:I\sn\times\sn(\sp\upsilon_s\varPi\sp)\to\upsilon_s\varPi$ for which $
x\ssp(\sp t\ssp,\cdot\sp)$ is holomorphic for all \œ$t\in I\sp$. We also
discuss aspects of the application of this result to establishing a
holomorphic solution map \œ$(\sp\xi\ssp,\varphi\sp)\mapsto y$ for functions \œ$
y:I\to\upsilon_s\varPi$ satisfying the ordinary differential equation $
y\ssp'=\varphi\circ[\,\roman{id\,},y\,]$ with initial condition $
y\ssp(\sp t_0)=\xi\,$.                                        \end{abstract}
%


In \cite{MH}\ssp, the following problem was considered. Fix a vector $\xi\ar 0
$ in a given complex Banach space $\varPi\sp$, and let $B\ar 0$ be an open
ball in $\varPi$ centered at $\xi\ar 0\ssp$. Consider differentiable curves $y
$ in $B\ar 0$ defined on the compact real interval \œ$I=[\,t\ar 0-\smb A\ssp,
t\ar 0+\smb A\,]$ and satisfying the differential equation \œ$\ssp
\roman E\ssp(\sp x\ssp,y\sp):y\ssp'(t)=x\ssp(\sp t\ssp,y\ssp(t))$ for all \œ$
t\in I$ with initial condition \œ$y\ssp(\sp t\ar 0)=\xi\ar 0\ssp$, where $x$
is a suitable function defined on \œ$I\snn\times\sn B\ar 0$ and having values
in $\varPi\sp$. Theorem 2 in \cite[p.\ 85\,]{MH} gave the result that if $
x\ar 0$ is suitably small in a certain Banach space $E\ar 1$ of functions $
x\ssp$, then an open neighbourhood $U$ of $x\ar 0$ in $E\ar 1$ exists such
that for all \œ$x\in U$ there is a unique $y$ with $\ssp
\roman E\ssp(\sp x\ssp,y\sp)\,$, and the function \œ$U\owns x\mapsto y$
defines a holomorphic map \œ$E\ar 1\to C\ssp(\spp I\sp,\varPi\sp)\,$. See
Constructions \ref{solu map} below and the discussion next to them for the
precise formulation.

It is essential for the preceding result that the map \œ$(\sp x\ssp,y\sp)
\mapsto x\circ[\,\sp\roman{id\,},y\,]:I\owns t\mapsto$ $x\ssp(\sp t\ssp,
y\ssp(t))$ be holomorphic. The purpose of this note is to establish this in a
setting more general than the one considered in \cite{MH}\ssp, see
Theorem \ref{map holom} below. We also indicate the main steps of the proof of
Theorem \ref{holom solu map} below generalizing \cite[Theorem 2\,]{MH}\ssp.

We use the conventions of \cite{Hi}\ssp. Therefrom, in particular, we recall
that $\biit R$ and $\biit C$ are the standard real and complex topological
fields, respectively. The topology of a topological vector space $E$ is $
\taurd E\ssp$, and its underlying set is $\vecs E\,$. Its filter of zero
neighbourhoods is $\ymp E\ssp$, and $\rajou E$ is the set of bounded sets. The
class of complex Banach\sp(able topological vector) spaces is $
\BaS(\biit C\ssp)\,$.\vskip.4mm

We also recall from \cite[Section 3\,]{Hi} that the particular holomorphy
class $\Cal H_{_T}$ has as its members exactly the maps \œ$\tilde f=(E\ssp,
F\spp,f\sp)$ of complex Hausdorff locally convex spaces $E\ssp,F$ where we
have $F$ locally (= Mackey, see \cite[p.\ 196\,]{Jr} or
\cite[Lemma 2.2, p.\ 15\,]{KM}\sp) complete, and $f$ a function with \œ$
\dom\sn f\in\taurd E\ssp$, i.e.\ open set in $E\ssp$, and\biggerlineskip2 \œ$
\rng\snn f\inc\vecs F\sp$, and $f$ continuous \œ$\taurd E\to\taurd F$ and $
\tilde f$ directionally differentiable, the last one meaning that for all
fixed \œ$x\in\dom\sn f$ and \œ$u\in\vecs E$ the limit $
\delta f\sp\fvalue(\sp x\ssp,u\sp)$ of $
t^{\sp\mminus 1}\sp(\sp f\sp\fvalue(\sp x+t\,u\sp)-f\sp\fvalue x\sp)$ as $
\Ce\setminus\ssn\{\sp 0\sp\}\owns t\to 0$ exists in the space $F\sp$.

Recall that we let $f\sp\fvalue x$ be the function value of $f$ at $x\ssp$,
and that \œ$f\,[\,A\,]=f\ssp\image\ssn A=$ \œ$\{\,y:\exi{x\in A}\,(\sp x\ssp,
y\sp)\in f\,\}\,$. In \cite{Hi}\ssp, we also agreed on the definitions \œ$
f\sp\inve=\{\ssp(\sp y\ssp,x\sp):$ $(\sp x\ssp,y\sp)\in f\,\}$ and $
\dom\sn f=\{\,x:\exi y\,(\sp x\ssp,y\sp)\in f\,\}$ and $\roman{dom}\yr 2 f=
\dom\sn(\dom\sn f\sp)\,$.

For \œ$\smb R\in\Rep=\{\,t:t\text{ real and }0<t\,\}\,$, here let \œ$
\bar{\roman D}\ssp(\smb R)$ be the set of \œ$t\in\Ce=\vecs\biit C$ with \œ$
|\ssp t\ssp|=\big((\sp\Rlp t\sp)^{\,2}+
                  (\sp\Imp t\sp)^{\,2}\big){\sp^{\frac 12}}\le\smb R\ssp$, and
put \œ$\roman D\ssp(\smb R)=\{\,t\in\Ce:|\ssp t\ssp|<\smb R\,\}\,$. Also let\biggerlineskip5
\œ$\roman S\snn\yr 1=\{\,\roman e^{\,2\ssp\pi\ssp\roman i\ssp t}\ssn:
t\in\Re\,\}\,$. If $F$ is a complex Hausdorff locally convex space, and if \œ$
\seq{\,\gamma\ssp(\zeta):\zeta\in\roman S\snn\yr 1\ssp}$ is a function defined
on $\roman S\snn\yr 1$ and continuous \œ$\taurd\spp\biit C\to\taurd F\sp$, we\biggerlineskip3
then let \œ$\oint\gamma\ssp(\zeta)\sp\,d\ssp\zeta=y$ if and only if \œ$y\in
\vecs F$ and for all \œ$\ell\in\Cal L\,(\spp F\sp,\biit C\ssp)\,$, i.e.\ for
all continuous linear functionals in $F\sp$, we have $\int_{\sp 0}^{\ssp 1}
\ell\ssp\fvalue(\sp\gamma\ssp(\sp\roman e^{\,2\ssp\pi\ssp\roman i\ssp s\sp}))
\,\roman e^{\,2\ssp\pi\ssp\roman i\ssp s}\,d\ssp s=\ell\ssp\fvalue v\ssp$. \vskip.5mm

Similarly, for \œ$k\in\N$ putting \œ$\roman S^{\ssp k}=
(\sp\roman S\snn\yr 1\snn)^{\,k}$, if \œ$f=\seq{\,\gamma\ssp(\bit\zzeta\sp):
\bit\zzeta\in\roman S^{\ssp k}\ssp}$ is a function\biggerlineskip4 defined on
$\roman S^{\ssp k}$ and continuous \œ$\taurd\spp\biit C^{\,k}\to\taurd F\sp$,
we let \œ$^k\!\oint\gamma\ssp(\bit\zzeta\sp)\sp\,d\ssp\bit\zzeta=y$ if and
only\biggerlineskip3 if \œ$y\in\vecs F$ and for all \œ$\ell\in
\Cal L\,(\spp F\sp,\biit C\ssp)$ we have \œ$
\int_{\sp[\sp 0\sp,\sp 1\sp]\sp^k}(\sp\roman m\sn\cdot\sn(\ssp\ell\circ f\sp))
\circ\roman p=\ell\ssp\fvalue y$ where\biggerlineskip8 $\roman m=
\Seq{\ssp\prod_{\,i\ssp\in\ssp k}\sp(\sp\bit\zzeta\fvalue i\sp):\bit\zzeta\in
\roman S^{\ssp k}\ssp}\,$ and $\,\roman p=\Seq{\ssp\seq{\,\exp\ssp(\sp
2\,\pi\,\roman i\,(\sp\bit\eeta\fvalue i\sp)):i\in k\,}:\bit\eeta\in\Re\,^k\ssp}\,$. \vskip.7mm

In \cite[Theorem 3.8\,]{Hi} we proved that \œ$\Cal H_{_T}\inc
\CinftyPi(\biit C\ssp)\,$. From this, assuming that \œ$(E\ssp,F\spp,f\sp)\in
\Cal H_{_T}\ssp$, it follows for every \œ$k\in\No$ that the $k^{\,\sixroman{th}}$
order {\it differential\ssp} $d^{\sp\,k\sn}f$ is defined on all of $\dom\sn f$
and has range included in the set $\Cal L^{\ssp k}(E\ssp,F\sp)$ of continuous
(complex-)\ssp multilinear maps \œ$E^{\,k}\to F\sp$. For fixed \œ$x\in
\dom\sn f$ and \œ$\biit u\in(\sp\vecs\varPi\sp)^{\,k}$ such that \œ$x+
\sum_{\,i\ssp\in\ssp k}\sp(\sp\bit\zzeta\fvalue i\sp)\,(\sp\biit u\fvalue i\sp)
\in\dom\sn f$ for all \œ$\bit\zzeta\in\bar{\roman D}\ssp(1)^{\,k}$, by
induction on \œ$k\in\No\ssp$, using composition by \œ$\ell\in
\Cal L\,(E\ssp,\biit C\ssp)\,$, the ordinary Cauchy\,--\,formula gives\vskip.7mm\centerline{$
\delta^{\,k\sn}f\sp\fvalue\big(\spp\seq{\ssp x\ssp}\concc\biit u\sp\big)=
d^{\sp\,k\sn}f\fvalue x\fvalue\biit u=
{^k}\!\oint\prod_{\,i\ssp\in\ssp k}\sp
(\sp\bit\zzeta\fvalue i\sp)^{\sp\mminus 2}\sp
f\sp\fvalue\big(\sp x+\sum_{\,i\ssp\in\ssp k}\sp
(\sp\bit\zzeta\fvalue i\sp)\,(\sp\biit u\fvalue i\sp)\big)
\sp\,d\ssp\bit\zzeta$}\vskip.7mm

\noin where $\delta^{\,k\sn}f$
      is the $k^{\,\sixroman{th}}$ order {\it variation\ssp} of $f\sp$.

Note that \hfill $[\,\sp\roman{id\,},y\,]=\{\ssp(\sp t\ssp,(\sp t\ssp,\xi\sp))
:(\sp t\ssp,\xi\sp)\in y\,\}\,$, \hfill and that we have for example\par$\mhyppy{4}
 \partial\ar 2\ssp x\circ[\,\sp\roman{id\,},y\,]\,.\,v\fvalue t
=(\sp\partial\ar 2\ssp x\circ[\,\sp\roman{id\,},y\,]\,.\,v\sp)\fvalue t
=(((\sp\partial\ar 2\ssp x\sp)\circ[\,\sp\roman{id\,},
                           y\,]\sp)\,.\,v\sp)\fvalue t$\par$\mhyppy{30.6}
=((\sp\partial\ar 2\ssp x\sp)\fvalue(\sp t\ssp,
      y\fvalue t\sp))\fvalue(\sp v\fvalue t\sp)
=(\sp\delta\ssp(\sp x\ssp(\sp t\ssp,\cdot\sp)))\fvalue(\sp y\fvalue t\ssp,
 v\fvalue t\sp)$ \hfill                and            \vskip.3mm\noin$
\partial_{\sixroman 2}^{\,2.}\sp x\circ[\,\sp\roman{id\,},y\,]\,.\,
                     [\,v\ssp,v\,]_{_{\roman{f1}}}\KN{1.3}\fvalue t
=(\sp\delta^{\,2.}\sp(\sp x\ssp(\sp t\ssp,\cdot\sp)))\fvalue\seq{\,
   y\fvalue t\ssp,v\fvalue t\ssp,v\fvalue t\,}\,$ under suitable conditions.

Above (and below) note that we make an explicit distinction between the domain
and range of the natural bijection \œ$\beta:\No\to\Zep\inc\Repp=\Rep\snn\cup
\{0\}\,$, and that\linebreak \œ$0\adot\ssp,1\adot\ssp,2\adot\in\No$ whereas \œ$
0\,,1\ssp,2\in\Zep\ssp$. For \œ$k\in\No$ we have \œ$k\ydot=\beta\fvalue k\in
\Zep$ and $k=(\sp k\ydot)\adot\ssp$. Further, for example $0\adot=\emptyset$
and $1\adot=\{\emptyset\}$ and $2\adot=\{\ssp\emptyset\,,1\adot\sp\}\,$, etc.

\begin{construction}\label{map}
With the data $I\sp,\varPi\sp,O\ssp$, we associate the map \œ$
(E\sp\sqcap F\sp,F\spp,f\sp)$ as follows. Here we assume that $I$ is a compact
real interval with nonempty interior, and that \œ$\varPi\in\BaS(\biit C\ssp)$
and \œ$O\in\Cal T\ssn\ar 0$ where \œ$\Cal T\ssn\ar 0$ is the trace on the set
\œ$I\snn\times\sn(\sp\vecs\varPi\sp)$ of the product topology of the
topologies $\taurd\spp\biit R$ and $\taurd\spp\varPi\sp$.

We put \œ$F=C\ssp(\spp I\sp,\varPi\sp)\,$, the complex Banach space of
continuous functions \œ$y:I\to$ $\vecs\varPi\sp$. Putting \œ$\Cal B=\{\,B:B
\inc O\text{ and }\rng B\in\rajou\varPi\text{ and }\exi{V\in\ymp\varPi}\,
B\spp\ai V\inc O\,\}\,$,\linebreak where we have \œ$B\spp\ai V=
\{\ssp(\sp t\ssp,\xi+\xi\ar 1):(\sp t\ssp,\xi\sp)\in B\text{ and }\xi\ar 1\in
V\,\}\,$, we let\vskip.5mm$\mhyppy{11}
S=\{\,x:x\in(\sp\vecs\varPi\sp)^{\,O}\text{ and }x$ continuous $
\Cal T\ssn\ar 0\to\taurd\spp\varPi$ and \newline\null\hfill$
\all{t\in I}\,(\varPi\sp,\varPi\sp,x\ssp(\sp t\ssp,\cdot))\in\Cal H_{_T}\text{
and }\sp\all{B\in\Cal B}\,x\image\snn B\in\rajou\varPi\,\}\,$.\KP{8}\linebreak
Writing \œ$N\sp(B\ssp,V\sp)=\{\,x\in S:x\image\snn B\inc V\,\}\,$, we let $E$
be the unique complex topological vector space (cf.\
\cite[p.\ 43\,--\,44\,]{Jr}\sp) with \œ$\seq{\,x\fvalue\smb P:x\in S\,}$
linear \œ$\sigrd E\to\sigrd\spp\varPi$ for all $\smb P\in O\ssp$, and such
that $\{\,N\sp(B\ssp,V\sp):B\in\Cal B\text{ and }V\in\ymp\varPi\,\}$ is a
filter basis for $\ymp E\,$. We finally put $\,f=\{\ssp(\sp x\ssp,y\ssp,z\sp):
x\in\vecs E\text{ and }y\ssp,z\in\vecs F$ $\text{and }z=
x\circ[\,\roman{id\,},y\,]\,\}\,$.                        \end{construction}

Holomorphy of the map \œ$\tilde f=(E\sp\sqcap F\sp,F\sp,f\sp)$ above is needed
when for fixed \œ$t\ar 0\in\Re$ and \œ$\smb A\in\Rep$ one wishes to establish
a holomorphic solution map \œ$(\sp\xi\ar 0\ssp,x\sp)\mapsto y$ for the
differential equations \œ$y\ssp'=x\circ[\,\sp\roman{id\,},y\,]$ with initial
condition \œ$y\fvalue t\ar 0=\xi\ar 0$ for $C^{\ssp 1.}$ functions \œ$y:I=
[\,t\ar 0-\smb A\ssp,t\ar 0+\smb A\,]\to\vecs\varPi\sp$, see
Theorem \ref{holom solu map} below. The space $E$ is Fr\'echet but as we have
no use of this, we don't take any pains for proving it.

\begin{theorem}\label{map holom}
In {\ssp\rm Construction \ref{map}\sp,} it holds that $
         (E\sp\sqcap F\sp,F\sp,f\sp)\in\Cal H_{_T}\ssp$.       \end{theorem}

\begin{proof} Noting that $E$ is Hausdorff locally convex, and that $F$ is
Banach, hence complete, so locally complete, we have to verify openness of
$\dom\sn f\sp$, and continuity and directional differentiability of $\tilde f\sp
$. To establish these, we first prepare the following tools. Fixing a
compatible norm \œ$\Nu:\vecs\varPi\owns\xi\mapsto|\ssp\xi\ssp|$ for $\varPi\sp
$, we let \œ$\ssp\roman B\ssp(\smb R)=\Nu\sp\inve\image[\,0\,,\smb R\sp\,{[}\ssp
$ and \œ$\ssp\bar{\roman B}\ssp(\smb R)=\Nu\sp\inve\image[\,0\,,\smb R\,]$ and
\œ$N\sn\ar 0\sp(B\ssp,\eps\sp)=N\ssp(B\ssp,\roman B\ssp(\eps))$ and\linebreak
\œ$N\sn\ar 1\spp(\delta\sp)=(\sp\vecs F\sp)\cap
(\sp\Pows(\spp I\snn\times\roman B\ssp(\delta\sp)))\,$. Note that \œ$
\{\ssp(\sp\vecs F\sp)\cap(\sp\Pows(\sp y\sbi{\,\roman B\sp(\eps)}\sp)):\eps\in
\Rep\ssp\}$ is a basis for the filter of $\taurd F\,$--\,neighbourhoods of any
\œ$y\in\vecs F\sp$, and that we agree on the definition $\sp
\Pows S=\{\,Q:Q\inc S\,\}\,$, power set.

Letting \œ$G=E\sp\sqcap F\sp$, the product topological vector space of $E$ and
$F\sp$, for the proof of \œ$\dom\sn f\in\taurd G\sp$, we first note that \œ$w
\in\dom\sn f$ if and only if there are $x\ssp,y\ssp,z$ with \œ$w=
(\sp x\ssp,y\sp)\in\vecs G$ and \œ$z=x\circ[\,\roman{id\,},y\,]\in\vecs F\sp$.
This implies that \œ$\dom\sn(\sp x\circ[\,\roman{id\,},y\,]\sp)=$ \œ$\dom z=I$
whence \œ$y=\rng[\,\sp\roman{id\,},y\,]\inc\dom x=O\ssp$. Conversely, if \œ$x
\in\vecs E$ and\linebreak \œ$y\in\vecs F$ with \œ$y\inc O\ssp$, for \œ$z=
x\circ[\,\roman{id\,},y\,]$ we then have \œ$\dom z=$ \œ$\dom y=I$ and $z$\linebreak
continuous \œ$\taurd\spp\biit R\to\taurd F\sp$, whence we get \œ$z\in\vecs F\sp
$. Consequently, we have \œ$\dom\sn f=$ $\{\ssp(\sp x\ssp,y\sp):x\in\vecs E\text{
and }y\in\vecs F\text{ and }y\inc O\,\}\,$.

By the above, to prove that \œ$\dom\sn f\in\taurd G\sp$, for arbitrarily fixed
\œ$y\in\vecs F$ with \œ$y\inc O\ssp$, it suffices to show existence of \œ$\eps
\in\Rep$ with \œ$y\sbi{\,\roman B\sp(\eps)}\inc O\ssp$. To obtain this, we\biggerlineskip3
note that as $y$ is continuous \œ$\taurd\spp\biit R\to\taurd F\sp$, so is \œ$
[\,\roman{id\,},y\,]:\taurd\spp\biit R\to\Cal T\ssn\ar 0\ssp$, and
consequently the set \œ$y=\rng[\,\roman{id\,},y\,]$ is $\Cal T\ssn\ar 0\ssp
$--\,compact.

Writing \œ$A\ssp(\sp t\ssp,\eps\sp)=\{\ssp(\sp t+s\ssp,y\fvalue t+\xi\sp):t+s
\in I\text{ and }|\ssp s\ssp|<\eps\text{ and }\xi\in\roman B\ssp(\eps)\ssp\}\,
$, we\linebreak consider the class \œ$\Cal I=\{\,(\sp\eps\ssp,t\ssp,
A\ssp(\sp t\ssp,\eps\sp)):\eps\in\Rep\text{ and }t\in I\text{ and }
A\ssp(\sp t\ssp,2\,\eps\sp)\inc O\,\}\,$.\linebreak Using \œ$O\in
\Cal T\ssn\ar 0\ssp$, we see that \œ$y\inc\bigcup\,(\rng\Cal I\sp)\,$, whence
by compactness of $y$ there is a\linebreak finite \œ$\Cal I\ar 0\inc\Cal I$
with \œ$y\inc\bigcup\,(\rng\Cal I\ar 0)\,$. For \œ$\eps=
\min\ssp(\sp\roman{dom}\yr 2\ssp\Cal I\ar 0)$ we then have \œ$\eps\in\Rep$,
and to prove that \œ$y\sbi{\,\roman B\sp(\eps)}\inc O\ssp$, arbitrarily fixing
\œ$t\in I$ and \œ$\xi\in\roman B\ssp(\eps)\,$, we should have \œ$(\sp t\ssp,
y\fvalue t+\xi\sp)\in O\ssp$. To see that this indeed is the case, by \œ$y\inc
\bigcup\,(\rng\Cal I\ar 0)$ we first see existence of some $\eps\ar 1\sp,
t\ar 1$ with \œ$(\sp\eps\ar 1\sp,t\ar 1\sp,A\ssp(\sp t\ar 1\sp,\eps\ar 1))\in
\Cal I\ar 0$ and \œ$(\sp t\ssp,y\fvalue t\sp)\in A\ssp(\sp t\ar 1\sp,
\eps\ar 1)\,$. Then we have \œ$|\,y\fvalue t-y\fvalue t\ar 1\ssp|<\eps\ar 1\ssp
$, whence \œ$|\,y\fvalue t+\xi-y\fvalue t\ar 1\ssp|<\eps+\eps\ar 1\le\eps\ar 1
+\eps\ar 1=2\,\eps\ar 1\ssp$, and consequently $(\sp t\ssp,y\fvalue t+\xi\sp)
\in A\ssp(\sp t\ar 1\sp,2\,\eps\ar 1)\inc O\ssp$.

To prove that $f$ is continuous \œ$\taurd\spp G\to\taurd F\sp$, we arbitrarily
fix \œ$(\sp x\ssp,y\ssp,z\sp)\in f$ and \œ$\eps\in\Rep$, and proceed to show
existence of \œ$W\sn\in\ymp G$ with the property that\linebreak \œ$
f\sp\fvalue(\sp x+u\ssp,y+v\sp)-z\inc I\snn\times\sn\roman B\ssp(\eps)$
whenever \œ$(\sp u\ssp,v\sp)\in W_{\sp}$. We now consider the class\vskip.2mm$\mhyppy{2}
\Cal I=\{\,(\sp\delta\ssp,t\ssp,A\ssp(\sp t\ssp,\delta\sp)):\delta\in\Rep$ and
$t\in I$ and $A\ssp(\sp t\ssp,3\,\delta\sp)\inc O$\par\hyppy{48mm}
and $x\,[\,A\ssp(\sp t\ssp,2\,\delta\sp)\,]\inc
x\fvalue(\sp t\ssp,y\fvalue t\sp)+\roman B\ssp(\sp\frac14\,\eps\sp)\ssp\}\,$.\vskip.4mm

\noin Now using also continuity of $x\ssp$, we see that \œ$y\inc
\bigcup\,(\rng\Cal I\sp)\,$, and since $y$ is $\Cal T\ssn\ar 0\ssp$--\,compact
and \œ$\rng\Cal I\inc\Cal T\ssn\ar 0\ssp$, there is some finite \œ$\Cal I\ar 0
\inc\Cal I$ with \œ$y\inc\bigcup\,(\rng\Cal I\ar 0)\,$. With\biggerlineskip3 \œ$
\delta=\min\ssp(\domm\Cal I\ar 0)$ then taking \œ$W=
N\sn\ar 0\sp(\sp y\sbi{\,\roman B\sp(\delta)}\ssp,\frac12\,\eps\sp)\sn\times\sn
N\sn\ar 1\spp(\delta\spp)\,$, we have \œ$W\in\ymp G$\biggerlineskip5 if \œ$
y\sbi{\,\roman B\sp(\delta)}\in\Cal B\ssp$. As \œ$
\rng y\sbi{\,\roman B\sp(\delta)}\inc\rng y+\roman B\ssp(\delta\spp)\in
\rajou\varPi$ and \œ$
(\sp y\sbi{\,\roman B\sp(\delta)}\sp)\sbi{\,\roman B\sp(\delta)}\inc
y\sbi{\,\roman B\sp(2\sp\delta)}\,$,\linebreak to get \œ$
y\sbi{\,\roman B\sp(\delta)}\in\Cal B\ssp$, it suffices that \œ$
y\sbi{\,\roman B\sp(2\sp\delta)}\inc O\ssp$. To obtain this, arbitrarily
fixing\linebreak \œ$t\in I$ and \œ$\xi\in\roman B\ssp(\sp 2\,\delta\sp)\,$, we
should have \œ$(\sp t\ssp,y\fvalue t+\xi\sp)\in O\ssp$. By \œ$y\inc
\bigcup\,(\rng\Cal I\ar 0)\,$, there\linebreak are some $\delta\ar 1$ and $
t\ar 1$ with \œ$(\sp t\ssp,y\fvalue t\sp)\in A\ssp(\sp t\ar 1\sp,\delta\ar 1)$
and \œ$(\sp\delta\ar 1\sp,t\ar 1\sp,A\ssp(\sp t\ar 1\sp,\delta\ar 1))\in
\Cal I\ar 0\ssp$. Then we\linebreak have \œ$
|\,y\fvalue t-y\fvalue t\ar 1\ssp|<\delta\ar 1\ssp$, and hence \œ$
|\,y\fvalue t+\xi-y\fvalue t\ar 1\ssp|<\delta\ar 1\snn+2\,\delta\le
\delta\ar 1\snn+2\,\delta\ar 1=3\,\delta\ar 1$\linebreak whence \œ$
(\sp t\ssp,y\fvalue t+\xi\sp)\in A\ssp(\sp t\ar 1\sp,3\,\delta\ar 1)\inc O\ssp
$. We also have \œ$(\sp x\ssp,y\sp)+W\inc\dom\sn f$ since for\linebreak $
(\sp u\ssp,v\sp)\in W$ we have $y+v\inc y\sbi{\,\roman B\sp(\delta)}\inc
y\sbi{\,\roman B\sp(2\sp\delta)}\inc O$ by the preceding.

Having now established \œ$W\in\ymp G$ and \œ$(\sp x\ssp,y\sp)+W\inc\dom\sn f\sp
$, we arbitrarily fix \œ$(\sp u\ssp,v\sp)\in W$ and put \œ$\Delta\sn\yr 1=
u\circ[\,\sp\roman{id\,},y+v\,]$ and \œ$\Delta\sn\yr 2=x\circ
[\,\sp\roman{id\,},y+v\,]-x\circ[\,\sp\roman{id\,},y\,]\,$.\linebreak Then \œ$
f\sp\fvalue(\sp x+u\ssp,y+v\sp)-z=(\sp x+u\sp)\circ[\,\sp\roman{id\,},y+v\,]-
x\circ[\,\sp\roman{id\,},y\,]=\Delta\sn\yr 1+\Delta\sn\yr 2$, and\linebreak
hence it suffices for arbitrarily fixed \œ$t\in I$ to establish \œ$
\Delta^\iota\fvalue t\in\roman B\ssp(\sp\frac12\,\eps\sp)$ for \œ$\iota=1\sp,
2\ssp$.\biggerlineskip4 For \œ$\iota=1\sp$, this is immediate by our arrangent.
For \œ$\iota=2\ssp$, by \œ$y\inc\bigcup\,(\rng\Cal I\ar 0)\,$, there (again)
are some $\delta\ar 1$ and $t\ar 1$ with \œ$(\sp t\ssp,y\fvalue t\sp)\in
A\ssp(\sp t\ar 1\sp,\delta\ar 1)$ and \œ$(\sp\delta\ar 1\sp,t\ar 1\sp,
A\ssp(\sp t\ar 1\sp,\delta\ar 1))\in\Cal I\ar 0\ssp$. Then we have \œ$
|\,y\fvalue t+v\fvalue t-y\fvalue t\ar 1\ssp|\le
|\,y\fvalue t-y\fvalue t\ar 1\ssp|+|\,v\fvalue t\,|<
\delta\ar 1\snn+\delta\le\delta\ar 1\snn+\delta\ar 1=2\,\delta\ar 1\ssp$,
whence we further get \œ$(\sp t\ssp,y\fvalue t+v\fvalue t\sp)\in
A\ssp(\sp t\ar 1\sp,2\,\delta\ar 1)\,$, and consequently\vskip.2mm\centerline{$
\Delta\sn\yr 2\fvalue t=x\fvalue(\sp t\ssp,y\fvalue t+v\fvalue t\sp)\pm
x\fvalue(\sp t\ar 1\sp,y\fvalue t\ar 1)-x\fvalue(\sp t\ssp,y\fvalue t\sp)\in
\roman B\ssp(\sp\frac14\,\eps\sp)-\roman B\ssp(\sp\frac14\,\eps\sp)\inc
\roman B\ssp(\sp\frac12\,\eps\sp)\,$.}\vskip.4mm

Preliminary to directional differentiability, we first prove that for \œ$
(\sp x\ssp,y\sp)\in\dom\sn f$ and \œ$v\in\vecs F$ and \œ$a=
\partial\ar 2\ssp x\circ[\,\sp\roman{id\,},y\,]\,$, we have \œ$a\,.\,v\in\vecs
F\sp$, i.e.\ that $a\,.\,v$ is continuous \œ$\taurd\spp\biit R\to
\taurd\spp\varPi\sp$. To establish this, we arbitrarily fix \œ$t\in I$ and \œ$
\eps\in\Rep$, and using \œ$(\sp x\ssp,y\sp)\in\dom\sn f\in\taurd G$ from above,
choose \œ$\delta\ar 0\in\Rep$ so that \œ$(\sp x\ssp,y+
\zeta\,\delta\ar 0\sp v\sp)\in\dom\sn f$ for all \œ$\zeta\in
\bar{\roman D}\ssp(1)\,$. As the function \œ$(\sp\tau\sp,\zeta\sp)\mapsto
x\fvalue(\sp\tau\sp,y\fvalue\tau+\zeta\,\delta\ar 0\sp(\sp v\fvalue\tau\sp))$
is continuous on the compact set \œ$I\snn\times\sn\bar{\roman D}\ssp(1)\,$,
there is \œ$\delta\in\Rep$ such that for all \œ$\zeta\in\bar{\roman D}\ssp(1)$
and for all real $s$ with $|\ssp s\ssp|<\delta$ and $t+s\in I$ we have\par\centerline{$
   x\fvalue(\sp t+s\ssp,y\fvalue(\sp t+s\sp)+\zeta\,\delta\ar 0\sp
                                          (\sp v\fvalue(\sp t+s\sp)))
-  x\fvalue(\sp t\ssp,y\fvalue t+\zeta\,\delta\ar 0\sp(\sp v\fvalue t\sp))
\in\roman B\ssp(\sp\delta\ar 0\ssp\eps\sp)\,$.}

\noin For $s\in\Re$ with $t+s\in I$ and $|\ssp s\ssp|<\delta\ssp$, we then have

$  a\,.\,v\fvalue(\sp t+s\sp)-a\,.\,v\fvalue t
=  \delta\ar 0\KN1^{\mminus 1}\sp(\sp a\,.\,(\sp\delta\ar 0\sp v\sp)\fvalue
   (\sp t+s\sp)-a\,.\,(\sp\delta\ar 0\sp v\sp)\fvalue t\sp)$\vskip.2mm\centerline{$
=  \delta\ar 0\KN1^{\mminus 1}\sp\oint\zeta^{\ssp\mminus 2}\sp(\sp
     x\circ[\,\sp\roman{id\,},y+\zeta\,\delta\ar 0\sp v\,]\fvalue(\sp t+s\sp)
   - x\circ[\,\sp\roman{id\,},y+\zeta\,\delta\ar 0\sp v\,]\fvalue
                                                 t\sp)\sp\,d\ssp\zeta$}\vskip.3mm$\mhyppy{8.7}
\in\delta\ar 0\KN1^{\mminus 1}\sp\roman B\ssp(\sp\delta\ar 0\ssp\eps\sp)
\inc\roman B\ssp(\eps)\,$.

For directional differentiability, given \œ$(\sp x\ssp,y\ssp,z\ssp)\in f$ and
\œ$(\sp u\ssp,v\ssp)\in\vecs G$ and \œ$\eps\in\Rep$, using the above
established continuity of $f\sp$, we first choose \œ$\delta\ar 1\in\Rep$ so
that for \œ$t\in\roman D\ssp(\sp\delta\ar 1)$ we have \œ$
 u\circ[\,\sp\roman{id\,},y+t\,v\,]
-u\circ[\,\sp\roman{id\,},y\,]\inc
I\snn\times\sn\roman B\ssp(\sp\frac12\,\eps\sp)\,$. Then we choose\linebreak
\œ$\delta\ar 0\in{]}\,0\,,\delta\ar 1\ssp]$ so that for \œ$S=\bigcup\ssp\{\,
y+t\,v+\zeta\,v+\zeta\ar 1\ssp v:t\ssp,\zeta\,,\zeta\ar 1\in
\bar{\roman D}\ssp(\sp\delta\ar 0)\ssp\}$ we have\linebreak \œ$S\inc O\ssp$,
this being possible by the above established result that \œ$\dom\sn f\in
\taurd G\sp$. As\linebreak $S$ is $\Cal T\ssn\ar 0\ssp$--\,compact and $x$ is
continuous, now $x\image S$ is $\taurd\spp\varPi\,$--\,compact, and in partic-
ular \œ$x\image S\in\rajou\varPi\sp$. Then letting $C$ be the closed convex
hull of $x\image S\ssp$, we also have $C\in\rajou\varPi\sp$. Hence, there is $
\delta\in{]}\,0\,,\delta\ar 0\ssp]$ with $
\delta\,\delta\ar 0\KN1^{\mminus 2}\ssp C\inc
\bar{\roman B}\ssp(\sp\frac12\,\eps\sp)\,$.

We now fix $t\in\Ce$ with $0<|\ssp t\ssp|<\delta\ssp$, and for \vskip.2mm\centerline{$
\varDelta=t^{\sp\mminus 1}\sp(\sp f\sp\fvalue(\sp x+t\,u\ssp,y+t\,v\sp)-z\sp)
-u\circ[\,\sp\roman{id\,},y\,]
-\partial\ar 2\ssp x\circ[\,\sp\roman{id\,},y\,]\,.\,v$}\vskip.6mm

\noin proceed to show that $\varDelta\inc
I\snn\times\sn\roman B\ssp(\eps)\,$. Since we have\vskip.4mm$\mhyppy{3.6}
\varDelta=u\circ[\,\sp\roman{id\,},y+t\,v\,]-u\circ[\,\sp\roman{id\,},y\,]$\par$\mhyppy{28}
+t^{\sp\mminus 1}\sp(\sp x\circ[\,\sp\roman{id\,},y+t\,v\,]
                        -x\circ[\,\sp\roman{id\,},y\,]\sp)
        -\partial\ar 2\ssp x\circ[\,\sp\roman{id\,},y\,]\,.\,v\,$,\vskip.7mm

\noin with $\,u\circ[\,\sp\roman{id\,},y+t\,v\,]-u\circ[\,\sp\roman{id\,},y\,]
\inc I\snn\times\sn\roman B\ssp(\sp\frac12\,\eps\sp)\,$, arbitrarily fixing $
t\ar 0\in I\sp$, it suffices\biggerlineskip4 to obtain $\,
(\sp t^{\sp\mminus 1}\sp(\sp x\circ[\,\sp\roman{id\,},y+t\,v\,]
                            -x\circ[\,\sp\roman{id\,},y\,]\sp)
-\partial\ar 2\ssp x\circ[\,\sp\roman{id\,},y\,]\,.\,v\sp)\fvalue t\ar 0\in
  \bar{\roman B}\ssp(\sp\frac12\,\eps\sp)\,$.\vskip.3mm

\noin To establish this by an appeal to the mean value theorem, considering
the curve\vskip.3mm\centerline{$
c:[\,0\,,1\,]\owns s\mapsto(\sp t^{\sp\mminus 1}\sp(\sp
 x\circ[\,\sp\roman{id\,},y+s\,t\,v\,]
-x\circ[\,\sp\roman{id\,},y\,]\sp)
-s\,\partial\ar 2\ssp x\circ[\,\sp\roman{id\,},y\,]\,.\,v\sp)\fvalue t\ar 0$}\vskip.6mm

\noin in the space $\varPi\sp$, it suffices that $\rng c\,'\inc
\bar{\roman B}\ssp(\sp\frac12\,\eps\sp)\,$. For arbitrarily fixed $s\in I\sp$,
we have\vskip.7mm\centerline{$
c\,'(s)=(\sp\partial\ar 2\ssp x\circ[\,\sp\roman{id\,},y+s\,t\,v\,]
-\partial\ar 2\ssp x\circ[\,\sp\roman{id\,},y\,]\sp)\,.\,v\fvalue t\ar 0\,$,}\vskip.5mm

\noin and to obtain \œ$c\,'(s)\in\bar{\roman B}\ssp(\sp\frac12\,\eps\sp)\,$,
we make another application of the mean value theorem by considering the curve\vskip.3mm\centerline{$
c\sp\ar 1\sn:[\,0\,,1\,]\owns s\ar 1\mapsto
(\sp\partial\ar 2\ssp x\circ[\,\sp\roman{id\,},y+s\ar 1\ssp s\,t\,v\,]
  -\partial\ar 2\ssp x\circ[\,\sp\roman{id\,},y\,]\sp)\,.\,v\fvalue t\ar 0\,$.}\vskip.6mm

\noin It then suffices to have $\rng c\ar 1\KN1'\inc
\bar{\roman B}\ssp(\sp\frac12\,\eps\sp)\,$. For arbitrarily fixed $s\ar 1\in I\sp
$, we have\vskip.6mm\centerline{$
c\ar 1\KN1'(\sp s\ar 1)=s\,t\,\partial_{\sixroman 2}^{\,2.}x\circ
[\,\sp\roman{id\,},y+s\ar 1\ssp s\,t\,v\,]\,.\,[\,v\ssp,
                       v\,]_{_{\roman{f1}}}\KN{1.3}\fvalue t\ar 0\,$,}\vskip.6mm

\noin whence we are done once \œ$\,t\,\partial_{\sixroman 2}^{\,2.}x\circ
[\,\sp\roman{id\,},y+s\ar 1\ssp s\,t\,v\,]\,.\,[\,v\ssp,
                      v\,]_{_{\roman{f1}}}\KN{1.3}\fvalue t\ar 0\in
        \bar{\roman B}\ssp(\sp\frac12\,\eps\sp)\,$ is verified. We have \ $
  t\,\partial_{\sixroman 2}^{\,2.}x\circ[\,\sp\roman{id\,},y+
  s\ar 1\ssp s\,t\,v\,]\,.\,[\,v\ssp,v\,]_{_{\roman{f1}}}\KN{1.3}\fvalue t\ar 0$\vskip.3mm$\mhyppy{7}
= t\,\delta_{\sixroman 0}^{\sp\mminus 2}\sp\partial_{\sixroman 2}^{\,2.}x\circ
  [\,\sp\roman{id\,},y+s\ar 1\ssp s\,t\,v\,]\,.\,[\,\delta\ar 0\ssp v\ssp,
  \delta\ar 0\ssp v\,]_{_{\roman{f1}}}\KN{1.3}\fvalue t\ar 0$\vskip.7mm$\mhyppy{7}
= t\,\delta_{\sixroman 0}^{\sp\mminus 2}\sp\oint\oint
  \zeta^{\ssp\mminus 2}\sp\zeta_{\sixroman 1}^{\ssp\mminus 2}\sp
  x\circ[\,\sp\roman{id\,},y+s\ar 1\ssp s\,t\,v+\zeta\,\delta\ar 0\ssp v
  +\zeta\ar 1\ssp\delta\ar 0\ssp v\,]\fvalue t\ar 0\,
    d\ssp\zeta\ar 1\ssp d\ssp\zeta$\vskip.7mm$\mhyppy{7.2}
\in t\,\delta\ar 0\KN1^{\mminus 2}\ssp C
=   t\,\delta^{\sp\mminus 1}\sp\delta\,\delta\ar 0\KN1^{\mminus 2}\ssp C
\inc\delta\,\delta\ar 0\KN1^{\mminus 2}\ssp C
\inc\bar{\roman B}\ssp(\sp\frac12\,\eps\sp)\,$.                  \end{proof}

Note that we could have shortened the preceding proof a little since the proof
of continuity of $f$ also anew gave openness of $\dom\sn f\sp$. However, in
the proof of openness we did not (essentially) need continuity of $x$ which
was strongly used in the proof of continuity. Also observe that our
arrangement of the use of the \q{compactness argument} made it possible to
manage without the axiom of choice contrary to the usual arrangement where
(say) one first (by use of AC) chooses a family \œ$t\mapsto\delta\sbi t$ with
$y\inc\bigcup\ssp\{\,S\ssp(\sp t\ssp,\delta\sbi t):t\in I\,\}$ and then using
compactness of $I$ establishes existence of a finite $J\inc I$ with $y\inc
\bigcup\ssp\{\,S\ssp(\sp t\ssp,\delta\sbi t):t\in J\,\}\,$.

\begin{constructions}\label{solu map}
Let $\varPi\sp,O\ssp,I\sp,E\ssp,F$ be as in Construction \ref{map}\sp, and for
some unique \œ$t\ar 0\in\Re$ and \œ$\smb A\in\Rep$ let \œ$I=
[\,t\ar 0-\smb A\ssp,t\ar 0+\smb A\,]\,$. Let \œ$F\ssn\ar 1=
C^{\ssp 1.}(\spp I\sp,\varPi\sp)$ and \œ$E\ar 2=$\linebreak \œ$\varPi\sp\sqcap
E\,$. Let $\Nu$ be as in the proof of Theorem \ref{map holom}\ssp. Fix \œ$
\xi\yr 0\in\vecs\varPi$ and \œ$\smb B\in\Rep$, and put \œ$B\ar 0=
\{\,\xi\yr 0+\xi:\xi\in\roman B\ssp(\smb B)\ssp\}\,$. Let $E\ar 1$ be the
Banach space of bounded continuous functions \œ$\varphi:I\sn\times\sn B\ar 0
\to\vecs\varPi$ for which \œ$\{\ssp(\varPi\sp,\varPi\sp,\varphi\ssp(\sp t\ssp,
\cdot\sp)):t\in I\,\}\inc\Cal H_{_T}\ssp$. Hence $\taurd E\ar 1$ is the
topology of uniform convergence on $I\snn\times\sn B\ar 0\ssp$. Further, let\vskip.2mm$\mhyppy{1.5}
\Sigma\ssp\ =\{\ssp(\sp\varphi\ssp,y\sp):\varphi\in\vecs E\ar 1\text{ and }y
\in\vecs F\text{ and }y\inc I\sn\times\sn B$\par\hyppy{51mm}
and $y\ssp'=\varphi\circ[\,\roman{id\,},y\,]\text{ and }y\fvalue t\ar 0=
\xi\yr 0\ssp\}$\hfill and\par$\mhyppy{1.5}
\varSigma\snn\ar 1=\{\ssp(\sp\xi\,,\varphi\ssp,y\sp):(\sp\xi\,,\varphi\sp)\in
\vecs E\ar 2\text{ and }y\in\vecs F\ssn\ar 1\text{ and }y\inc O$\par\hyppy{51mm}
and $y\ssp'=\varphi\circ[\,\sp\roman{id\,},y\,]\text{ and }y\fvalue t\ar 0=
\xi\,\}$ \hfill and\par$\mhyppy{1.5}
f\sn\ar 1\sp\,=\big\{\ssp(\sp\xi\,,\varphi\ssp,y\ssp,z\sp):(\sp\xi\,,
\varphi\sp)\in\vecs E\ar 2\text{ and }y\ssp,z\in\vecs F\ssn\ar 1$\par\hyppy{54.7mm}
and $z=y-I\sn\times\sn\{\sp\xi\sp\}-\int_{\sp t_0}\varphi\circ
[\,\roman{id\,},y\,]\,\big\}\,$.\vskip.4mm

\noin Here we understand that $\int_{\sp t_0}v$ is defined on the largest
interval $J$ with \œ$t\ar 0\in J\inc$ $\dom v$ to the points of which $v$ can
be integrated. We also let \œ$\ell\ar 0=\Seq{\sp\int_{\sp t_0}v:v\in\vecs F\,}
$ and let \ssp P$\ssp(x)$ for \œ$x=(\sp\xi\,,\varphi\sp)$ mean that for some \œ$
\smb B\snn\ar 1\in\Rep$ and for \œ$A=\{\ssp(\sp t\ssp,\xi+\xi\ar 1):$ $t\in I\text{
and }\xi\ar 1\snn\in\roman B\ssp(\smb B\snn\ar 1)\ssp\}\,$, we have $A\inc O$
and $\sup\ssp(\sp\Nu\circ\varphi\image\ssn A\sp)
<\frac12\,\smb A^{\mminus 1}\sp\smb B\snn\ar 1\ssp$.     \end{constructions}

About the setting of Constructions \ref{solu map} above
\cite[Theorem 2, p.\ 85\,]{MH} asserts that if \œ$\varphi\ar 0\in\vecs E\ar 1$
with \œ$\sup\ssp(\rng\snn(\sp\Nu\circ\varphi\sp))<
\frac12\,\smb A^{\mminus 1}\sp\smb B\ssp$, there is $U$ with \œ$\varphi\ar 0
\in U\in\taurd E\ar 1$\biggerlineskip3 and \œ$U\inc\dom\Sigma$ and \œ$
(E\ar 1\sp,F\spp,\Sigma\,|\,U\sp)\in\Cal H_{_T}\ssp$. We note that without
introducing any essential complications in the required proof, as \œ$
\rng\Sigma\inc\vecs F\ssn\ar 1$ and \œ$\ell\ar 0\in\Cal L\,(\spp F\sp,
F\ssn\ar 1)\,$, this result even persists if we take $F\ssn\ar 1$ in place of
$F\sp$. We may further (formally) strengthen it by introducing the
reformulation (m) that \ $(E\ar 1\sp,F\ssn\ar 1\sp,\Sigma\sp)\in\Cal H_{_T}$
and\par\centerline{$
\{\,\varphi:\varphi\in\vecs E\ar 1\text{ and }\sup\ssp(\rng\snn(\sp\Nu\circ
\varphi\sp))<\frac12\,\smb A^{\mminus 1}\sp\smb B\,\}\inc\dom\Sigma\,$.}

Notice that the preceding inclusion \œ$\ldots\inc\dom\Sigma$ is obtained by
estimating the Lipschitz\,--\,constant for the proof of existence of the
solution to \œ$y\ssp'=\varphi\circ[\,\roman{id\,},y\,]$ and $y\fvalue t\ar 0=
\xi\yr 0\,$ from \ \ $\int_{t_0}(\sp\varphi\circ[\,\roman{id\,},u\,]-
               \varphi\circ[\,\roman{id\,},v\,]\sp)$\vskip.7mm\centerline{$
=\int_{t_0}\int_0^{\ssp 1}\partial\ar 2\sp\varphi\circ[\,\roman{id\,},
                   s\,u+(\spp 1-s\sp)\,v\,]\,.\,(\sp u-v\sp)\sp\,d\ssp s\,$,}\vskip.7mm\noin
and \ $
 \partial\ar 2\sp\varphi\circ[\,\roman{id\,},s\,u+(\spp 1-s\sp)\,v\,]\,.\,z
=\eps^{\sp\mminus 1}\oint\zeta^{\ssp\mminus 2}\ssp\varphi\circ[\,\roman{id\,},
 s\,u+(\spp 1-s\sp)\,v+\eps\,\zeta\,z\,]\sp\,d\,\zeta\,$,\par\noin
and observing that the solution exists on $I$ if with \œ$\smb M=
\sup\ssp(\rng\snn(\sp\Nu\circ\varphi\sp))$ we can choose \œ$\eps\ssp,\smb R\in
\Rep$ so that \œ$\smb R+\eps<\smb B$ and \œ$\smb A\,\smb M\le\smb R$ and \œ$
\eps^{\sp\mminus 1}\smb A\,\smb M<1\ssp$, which is possible exactly when $
\smb M<\frac12\,\smb A^{\mminus 1}\sp\smb B\ssp$.

Using the implicit function theorem in \cite[p.\ 235\,]{Hi-Stud M}\ssp,
restated in \cite[Theorem 4.1\ssp]{Hi}\ssp, we can improve (m) above to obtain
the following

\begin{theorem}\label{holom solu map}
In the setting of {\ssp\rm Constructions \ref{solu map}} above$\ssp,$ it holds
that\vskip.1mm\centerline{$
(E\ar 2\ssp,F\ssn\ar 1\sp,\varSigma\snn\ar 1)\in\Cal H_{_T}$ and $\sp
\{\,x:x\in\vecs E\text{ and P}\ssp(x)\ssp\}\inc\dom\varSigma\snn\ar 1\ssp$.}
                                                               \end{theorem}

\begin{proof} Since the arguments required have much overlap both with one
another and those already given in detail in the proof Theorem \ref{map holom}
above, as well as with those in \cite[Section 5\,]{Hi}\ssp, we only indicate
the main aspects of them. First, we note that to get \œ$
(E\ar 2\ssp,F\ssn\ar 1\sp,\varSigma\snn\ar 1)\in\Cal H_{_T}\ssp$, arbitrarily
fixing \œ$w=(\sp x\ssp,y\sp)=(\sp\xi\,,\varphi\ssp,y\sp)\in\varSigma\snn\ar 1\ssp
$, {\it if\linebreak $\varSigma\snn\ar 1$ is a function\ssp}, and as $
F\ssn\ar 1$ is Banach, hence locally complete, it suffices to show existence
of $W$ with \œ$w\in W$ and \œ$(E\ar 2\ssp,F\ssn\ar 1\sp,\varSigma\snn\ar 1\cap
W\sp)\in\CPi{1.}(\biit C\ssp)\,$. For this, as we have \œ$\varSigma\snn\ar 1=
f\snn\ar 1\KN{.7}\inve\image\snn\{\ssp\bnull F\}\,$, by
\cite[Theorem 4.1\ssp]{Hi} it suffices that \œ$
(E\ar 2\spp\sqcap F\ssn\ar 1\sp,F\ssn\ar 1\sp,f\snn\ar 1)\in
\CPi{1.}(\biit C\ssp)$ and that $\partial\ar 2\sp f\snn\ar 1\KN{.7}\fvalue w$
is bijective $\vecs F\ssn\ar 1\to\vecs F\ssn\ar 1\ssp$.

As for $\varSigma\snn\ar 1$ being a function, we note that if for \œ$\iota=
1\ssp,2$ we have \œ$(\sp\xi\,,\varphi\ssp,y\sbi\iota)\in\varSigma\snn\ar 1\ssp
$, by standard ordinary differential equation theory in Banach spaces, to get
\œ$y\ar 1=y\ar 2\ssp$, it suffices that for all \œ$\sp\smb P\sp=(\sp t\ssp,
\xi\ar 0)\in O$ there are some $N$ and \œ$\smb M\in\Rep$ with \œ$\smb P\in N
\in\Cal T\ssn\ar 0$ and also such that for all \œ$
(\sp t\ar 1\sp,\xi\ar 1)\ssp,(\sp t\ar 1\sp,\xi\ar 2)\in N$ we have \œ$
|\,\varphi\fvalue(\sp t\ar 1\sp,\xi\ar 1)
  -\varphi\fvalue(\sp t\ar 1\sp,\xi\ar 2)\,|\le$ \œ$
\smb M\,|\,\xi\ar 1\snn-\xi\ar 2\ssp|\,$. To establish this by way of applying
the mean value theorem, it suffices to observe that for small $\xi\ar 3\ssp$,
say for \œ$|\,\xi\ar 3\ssp|\le\delta\ar  0\ssp$, and for $(\sp t\ar 1\sp,
\xi\ar 4)$ close to $\smb P\,$, by the Cauchy\,--\,formula we have \œ$
\partial\ar 2\sp\varphi\fvalue(\sp t\ar 1\sp,\xi\ar 4)\fvalue\xi\ar 3=
\oint\zeta^{\ssp\mminus 2}\sp\varphi\fvalue(\sp t\ar 1\sp,
\xi\ar 4+\zeta\,\xi\ar 3)\,\sp d\ssp\zeta\,$, and that by continuity of $
\varphi\ssp$, and by local boundedness of $\varPi\sp$, the norm of this
remains bounded by some \œ$\smb M\snn\ar 0<\plusinfty$ if $(\sp t\ar 1\sp,
\xi\ar 4)$ is confined to a suitably small neighbourhood of $\smb P\,$. By
homogeneity of $\partial\ar 2\sp\varphi\fvalue(\sp t\ar 1\sp,\xi\ar 4)\,$,
taking \œ$\smb M=\delta\ar 0\KN1^{\mminus 1}\sp\smb M\snn\ar 0\ssp$, we then
have \œ$|\,\sp\partial\ar 2\sp\varphi\fvalue(\sp t\ar 1\sp,\xi\ar 4)\fvalue
\xi\ar 3\ssp|\le$ $
\delta\ar 0\KN1^{\mminus 1}\sp\smb M\snn\ar 0\ssp|\,\xi\ar 3\ssp|
=\smb M\,|\,\xi\ar 3\ssp|$ for arbitrary $\xi\ar 3\in\vecs\varPi\sp$.

As for \œ$(E\ar 2\spp\sqcap F\ssn\ar 1\sp,F\ssn\ar 1\sp,f\snn\ar 1)\in
\CPi{1.}(\biit C\ssp)\,$, note that introducing the continuous linear map \œ$
\ssp\roman s:F\ssn\ar 1\sqcap F\ssn\ar 1\sqcap\varPi\to F\ssn\ar 1$ defined by
\œ$(\sp v\ssp,v\ar 1\sp,\xi\ar 1)\mapsto
v-v\ar 1-I\snn\times\sn\{\ssp\xi\ar 1\}\,$, we have \œ$f\snn\ar 1\inc$ \œ$
\roman s\circ[\,\sp\roman{pr}\ar 2\ssp,\spp\ell\ar 0\circ f\circ
[\,\sp\roman{pr}\ar 2\circ\roman{pr}\ar 1\sp,\roman{pr}\ar 2\ssp]\ssp,
\roman{pr}\ar 1\snn\circ\roman{pr}\ar 1\ssp]\,$. Since also \œ$(\spp F\sp,
F\ssn\ar 1\sp,\spp\ell\ar 0)$ is a continuous linear map, and as by
Theorem \ref{map holom} above and by \cite[Theorem 3.8\,]{Hi} we have \œ$
(E\ssp,F\spp,f\sp)\in$ \œ$\Cal H_{_T}\inc\CinftyPi(\biit C\ssp)\,$, the chain
rule \cite[Proposition 0.11, p.\ 240\,]{Hi-Stud M} gives us the required
result that $(E\ar 2\spp\sqcap F\ssn\ar 1\sp,F\ssn\ar 1\sp,f\snn\ar 1)\in
\CinftyPi(\biit C\ssp)\inc\CPi{1.}(\biit C\ssp)\,$.

For bijectivity of \œ$\partial\ar 2\sp f\snn\ar 1\KN{.7}\fvalue w:\vecs
F\ssn\ar 1\to\vecs F\ssn\ar 1\ssp$, with \œ$a=\partial\ar 2\sp\varphi\circ
[\,\sp\roman{id\,},y\,]\,$, from the proof of Theorem \ref{map holom} above,
and from \œ$(E\sp\sqcap F\sp,F\spp,f\sp)\in\Cal H_{_T}\inc
\CinftyPi(\biit C\ssp)\inc\CPi{1.}(\biit C\ssp)\,$, we first obtain \œ$
\seq{\,a\,.\,v:v\in\vecs F\,}\in\Cal L\,(\spp F\sp,F\sp)\,$, and then \œ$
\partial\ar 2\sp f\snn\ar 1\KN{.7}\fvalue w=
\seq{\,v-\ell\ar 0\KN1\fvalue(\sp a\,.\,v\sp):v\in\vecs F\ssn\ar 1\ssp}\,$.
Hence we are reduced to proving that for all \œ$v\ar 1\in\vecs F\ssn\ar 1\ssp
$, in the space $\varPi$ the ordinary differential equation \œ$v\ssp'-a\,.\,v=
v\ar 1\KN1'$ with initial condition \œ$v\fvalue t\ar 0=
v\ar 1\KN1\fvalue t\ar 0$ has a unique solution $v\ssp$. Putting \œ$\hat a=
\{\ssp(\sp t\ssp,\xi\ar 1\sp,a\fvalue t\fvalue\xi\ar 1):t\in I\text{ and }\xi
\ar 1\in\vecs\varPi\,\}\,$, again by standard theory, it suffices that
$\hat a$ is continuous \œ$\Cal T\ssn\ar 0\to\taurd\spp\varPi$ and $
\bigcup\,(\rng a\sp)\image(\sp\roman B\ssp(1))\in\rajou\varPi\sp$.

For continuity of $\hat a$ note that when \œ$\,(\sp t\ssp,\xi\ar 0)\ssp,
(\sp s\ssp,\xi\ar 1)\in\Re\snn\times\sn(\sp\vecs\varPi\sp)$ with \œ$t\ssp,t+s
\in I\sp$, we have \œ$\,\hat a\fvalue(\sp t+s\ssp,\xi\ar 0+\xi\ar 1)-
\hat a\fvalue(\sp t\ssp,\xi\ar 0)=(\sp a\fvalue(\sp t+s\sp)-
(\sp a\fvalue t\sp))\fvalue\xi\ar 0+a\fvalue(\sp t+s\sp)\fvalue\xi\ar 1\ssp$.
As $\varphi$ is continuous at the points of the compact set $y\ssp$, some
small \œ$\delta\ar 0\in\Rep$ exists with \œ$\{\,\varphi\fvalue(\sp t\ar 1\sp,
y\fvalue t\ar 1\snn+\zeta\,\xi\ar 1):t\ar 1\in I\text{ and }\zeta\in
\roman S\snn\yr 1\text{ and }\xi\ar 1\in\roman B\ssp(\sp\delta\ar 0)\ssp\}\in
\rajou\varPi\sp$. Then using \œ$a\fvalue t\ar 1\KN1\fvalue\xi\ar 1=
\oint\zeta^{\ssp\mminus 2}\sp\varphi\fvalue(\sp t\ar 1\sp,y\fvalue t\ar 1\snn+
\zeta\,\xi\ar 1)\,\sp d\ssp\zeta\,$, we deduce both \œ$
\bigcup\,(\rng a\sp)\image(\sp\roman B\ssp(1))\in\rajou\varPi\sp$, and that
for the fixed $t$ we have \œ$|\,a\fvalue(\sp t+s\sp)\fvalue\xi\ar 1\sp|\le
\smb M\subtext{some}\ssp|\,\xi\ar 1\sp|\,$. To complete the proof of
continuity of $\hat a\,$, upon putting \œ$v=I\snn\times\sn\{\ssp\xi\ar 0\}\,$,
we have \œ$(\sp a\fvalue(\sp t+s\sp)-(\sp a\fvalue t\sp))\fvalue\xi\ar 0=$ $
a\,.\,v\fvalue(\sp t+s\sp)-a\,.\,v\fvalue t$ which becomes small with $s$
since we have $a\,.\,v\in\vecs F\sp$.                            \end{proof}

Even in the case \œ$O=I\snn\times\sn B\ar 0\ssp$, for nontrivial spaces $
\varPi\sp$, if we fix \œ$\xi=\xi\yr 0$, the asser- tion of
Theorem \ref{holom solu map} is strictly stronger than (m) above as, by
considering a sequence \œ$\seq{\,\{\ssp(\sp t\ssp,\xi\,,(\ssp\ell\ssp\fvalue
\xi\sp)^{\,i}\sp\xi\ar 1):(\sp t\ssp,\xi\sp)\in O\,\}:i\in\No\ssp}$ with \œ$
\xi\ar 1\in\vecs\varPi\setminus\sn\{\sp\bnull\varPii\}$ and \œ$\ell\in
\Cal L\,(\varPi\sp,\biit C\ssp)$\linebreak suitably chosen, the topology $
\taurd E\ar 1$ of uniform convergence on $O$ is strictly stronger than the
topology of uniform convergence on all sets \œ$B\in\Cal B$ of
Construction \ref{map} above. From the next example, by taking \œ$B\ar 0=
\roman B\ssp(1)\,$, it even follows that we have the strict inclusion $\vecs
E\ar 1\subset\vecs E$ in the case where $\varPi=
\ell\KP{.7}\RHB{.3}{^p\spp}(\No\ssp,\biit C\ssp)$ with $1\le p<\plusinfty\,$.

\begin{example}\label{entire fu}
We construct an entire function $\chi$ on a complex Banach space $\varPi$
which is not bounded in the open unit ball. Namely, assuming that \œ$1\le p<
\plusinfty\,$, we show this for \œ$\varPi=\ell\KP{.7}\RHB{.3}{^p\spp}(\No\ssp,
\biit C\ssp)$ and \œ$\chi=\Seq{\sp\sum_{\,i\ssp\in\ssp I\!\!N_{\roman o}}i\,(\sp
x\fvalue i\sp)^{\,i}\sn:x\in\vecs\varPi\,}\,$. Let \œ$\Nu=$ \œ$
\Seq{\,\big(\sp\sum_{\,i\ssp\in\ssp I\!\!N_{\roman o}}
|\,x\fvalue i\,|^{\,p}\big){\sp^{p^{-1}}}\!:x\in\vecs\varPi\,}\,$, and for
convenience write \œ$\|\ssp x\ssp\|=\Nu\sp\fvalue x\ssp$. \hfill
We have \ \
$
\sup\ssp\{\,|\,\chi\fvalue x\,|:x\in\Nu\sp\inve\image[\,0\,,1\,{[\sp}\,\}\ge
\sup\ssp\{\,i\,t^{\,i}\sn:i\in\N\text{ and }0<t<1\ssp\}=\plusinfty\,$.

To simplify the notations a bit, we agree to let \œ$x\sbi i=x\fvalue i$ and \œ$
x_i^{\,k}=(\sp x\fvalue i\sp)^{\,k}$. To see that \œ$\dom\chi=\vecs\varPi\sp$,
it suffices for \œ$a\in\Repp\KN1\potNo$ with
\œ$\sum_{\,i\ssp\in\ssp I\!\!N_{\roman o}}a_i^{\,p}<\plusinfty$ that we have
\œ$\sum_{\,i\ssp\in\ssp I\!\!N_{\roman o}}i\,a_i^{\,i}<\plusinfty\,$. Note
that there is some \œ$\smb N\in\No$ such that we have \œ$a\sbi i\le\frac12$
whenever $\smb N\inc i\in\No\ssp$. For $\smb M=2+\|\ssp a\ssp\|\,$, then\vskip.7mm\centerline{$
   \sum_{\,i\ssp\in\ssp I\!\!N_{\roman o}}i\,a_i^{\,i}
=  \sum_{\,i\ssp\in\sp\ssmb N}\sp i\,a_i^{\,i}
   +\sum_{\,i\ssp=\sp\ssmb N}^{\ \infty}i\,a_i^{\,i}
\le(\ssp\smb N\ydot-1\ssp)\,\sum_{\,i\ssp\in\sp\ssmb N}\smb M^{\,i}
   +\sum_{\,i\ssp=\sp\ssmb N}^{\ \infty}i\,2^{\sp\mminus i}$}\vskip.7mm$\mhyppy{15}
=  (\ssp\smb N\ydot-1\ssp)\,(\spp\smb M-1\spp)^{\sp\mminus 1}\sp
   (\spp\smb M^{\,\ssmb N}\snn-1\spp)
   +\frac 12\,\seq{\,t^{\,\ssmb N}\sp(\spp 1-t\sp)^{\sp\mminus 1}\sn:
                   t\in{]}\,0\,,1\,{[}\ssp\,}\ssp'(\frac 12)$\vskip.7mm$\mhyppy{15}
=  (\ssp\smb N\ydot-1\ssp)\,(\spp\smb M-1\spp)^{\sp\mminus 1}\sp
   (\spp\smb M^{\,\ssmb N}\snn-1\spp)
   +2^{\sp\,1\sp-\sp\ssmb N\yydot}(\spp 1+\smb N\ydot\sp)<\plusinfty\,$.\vskip.7mm

To prove that $\chi$ is continuous \œ$\taurd\spp\varPi\to\taurd\spp\biit C$ at
an arbitrarily fixed \œ$x\in\vecs\varPi\sp$, putting \œ$\smb M=
1+\|\ssp x\ssp\|\,$, and considering \œ$u\in\vecs\varPi$ with \œ$
\|\ssp u\ssp\|<\frac 18\,$, there is \œ$\smb N\in\No$ such that for \œ$\smb N
\inc i\in\No$ we have \œ$|\,x\sbi i\ssp|<\frac 18\,$. Noting that \œ$
\sup\ssp\{\,i\,2^{\sp\mminus i}\sn:i\in\No\ssp\}\le$\biggerlineskip4 $
(\sp\ln 4\sp)^{\sp\mminus 1}<1\ssp$, now \
$\chi\fvalue(\sp x+u\sp)-\chi\fvalue x
= \sum_{\,i\ssp\in\ssp I\!\!N_{\roman o}}i\,((\sp x\sbi i+u\sbi i\spp)^{\,i}-
   x_i^{\,i}\sp)$\vskip.7mm$\mhyppy{20}
= \sum_{\,i\ssp\in\sp\ssmb N}\sp i\,((\sp x\sbi i+u\sbi i\spp)^{\,i}-
   x_i^{\,i}\sp)+\sum_{\,i\ssp=\sp\ssmb N}^{\ \infty}i\,((\sp x\sbi i+
   u\sbi i\spp)^{\,i}-x_i^{\,i}\sp)\,$,\hfill whence \vskip.7mm$\mhyppy{3}
   |\,\chi\fvalue(\sp x+u\sp)-\chi\fvalue x\,|
\le\sum_{\,i\ssp\in\sp\ssmb N}i\,\|\ssp u\ssp\|\,(\sp i\,\smb M^{\,i\yydot-\sp 1}\spp)
 + \sum_{\,i\ssp=\sp\ssmb N}^{\ \infty}
 i\,\|\ssp u\ssp\|\,(\sp i\,4^{\,\sp 1\sp-\ssp i\yydot}\sp)$\vskip.7mm$\mhyppy{32}
=  \|\ssp u\ssp\|\ssp\big(\spp\sum_{\,i\ssp\in\sp\ssmb N}
    i^{\,2.}\sp\smb M^{\,i\yydot-\sp 1}
  + \sum_{\,i\ssp=\sp\ssmb N}^{\ \infty}
   i^{\,2.}\,4^{\,\sp 1\sp-\ssp i\yydot}\sp\big)$\vskip.7mm$\mhyppy{32}
\le\|\ssp u\ssp\|\ssp\big(\sp\smb N^{\,3.}\sp\smb M^{\,\ssmb N\yydot-\sp 2}
  + 2\,\sum_{\,i\ssp=\sp\ssmb N}^{\ \infty}
    i\,2^{\,\sp 1\sp-\ssp i\yydot}\sp\big)$\vskip.7mm$\mhyppy{32}
=  \|\ssp u\ssp\|\,(\sp\smb N^{\,3.}\sp\smb M^{\,\ssmb N\yydot-\sp 2}
  + 2\,\seq{\,t^{\,\ssmb N}\sp(\spp 1-t\sp)^{\sp\mminus 1}\sn:t\in
    {]}\,0\,,1\,{[}\ssp\,}\ssp'(\frac 12))$\vskip.5mm$\mhyppy{32}
=  \|\ssp u\ssp\|\,(\sp\smb N^{\,3.}\sp\smb M^{\,\ssmb N\yydot-\sp 2}
  + 2^{\sp\,3\ssp-\sp\ssmb N\yydot}(\spp 1+\smb N\ydot\spp))\,$,\vskip.5mm

\noin and continuity at $x$ then follows.

To prove that $(\varPi\sp,\biit C\ssp,\chi\sp)$ is directionally
differentiable with variation\par\centerline{$
\delta\ssp\chi=\Seq{\,\sum_{\,i\ssp\in\ssp I\!\!N_{\roman o}}
i^{\,2.}(\sp\sigrd z\fvalue i\sp)^{\,i\yydot-\sp 1}\sp(\sp\taurd z\fvalue
i\sp):z\in\vecs(\varPi\spp\sqcap\varPi\sp)\,}\,$,}\vskip.3mm

\noin where note that for example \œ$\sigrd(\sp x\ssp,u\sp)\fvalue i=x\sbi i\ssp
$, arbitrarily fixing \œ$x\ssp,u\in\vecs\varPi$ and\biggerlineskip3 $\eps\in
\Rep$, and writing \hfill $
\varDelta\ssp(t)=\sum_{\,i\ssp\in\ssp I\!\!N_{\roman o}}
i\,|\,t^{\sp\mminus 1}\sp((\sp x\sbi i+t\,u\sbi i\spp)^{\,i}-
x_i^{\,i}\sp)-i\,x_i^{\,i\yydot-\sp 1}\sp u\sbi i\,|\,$,\hfill\null\vskip.5mm

\noin it suffices to show existence of \œ$\delta\in{]}\,0\,,1\,]$ such that
for all \œ$t\in\Ce$ with \œ$0<|\ssp t\ssp|<\delta$ we have \œ$\varDelta\ssp(t)
<\eps\,$. To see this, let \œ$\smb M=1+\|\ssp x\ssp\|+\|\ssp u\ssp\|\,$, and
again note that there is $\smb N\in\No$ such that for $\smb N\inc i\in\No$ we
have $|\,x\sbi i\ssp|\ssp,|\,u\sbi i\ssp|<\frac 18\,$. Then we may compute\vskip.7mm$\mhyppy{10}
\varDelta\ssp(t)=\sum_{\,i\ssp\in\ssp I\!\!N_{\roman o}}
i^{\,2.}\sp\big|\,\int_{\sp 0}^{\ssp 1}
((\sp x\sbi i+s\,t\,u\sbi i)^{\,i\yydot-\sp 1}\snn-
x_i^{\,i\yydot-\sp 1}\sp)\,u\sbi i\sp\,d\ssp s\,\big|$\vskip.7mm$\mhyppy{17.5}
=  \sum_{\,i\ssp\in\ssp I\!\!N_{\roman o}}
   i^{\,2.}\sp(\sp i\ydot\snn-1\spp)\,\big|\,t\,u_i^{\,2}\ssp
\int_{\sp 0}^{\ssp 1}\int_{\sp 0}^{\ssp 1}s\,
(\sp x\sbi i+s\ar 1\ssp s\,t\,u\sbi i)^{\,i\yydot-\ssp 2}
\,d\ssp s\ar 1\ssp d\ssp s\,\big|$\vskip.7mm$\mhyppy{17.5}
\le|\ssp t\ssp|\,\|\ssp u\ssp\|^{\,2}\ssp
   \sum_{\,i\ssp=\ssp 2.}^{\ \infty}
i^{\,3.}\ssp
\int_{\sp 0}^{\ssp 1}\int_{\sp 0}^{\ssp 1}
|\,x\sbi i+s\ar 1\ssp s\,t\,u\sbi i\ssp|^{\,i\yydot-\ssp 2}
\,d\ssp s\ar 1\ssp d\ssp s$\vskip.7mm$\mhyppy{17.5}
\le|\ssp t\ssp|\,\|\ssp u\ssp\|^{\,2}\ssp
   \sum_{\,i\ssp=\ssp 2.}^{\ \infty}
i^{\,3.}\ssp\big(\spp
|\,x\sbi i\ssp|+|\,u\sbi i\ssp|\spp\big)
{^{\,i\yydot-\ssp 2}}$\vskip.7mm$\mhyppy{17.5}
\le|\ssp t\ssp|\,\|\ssp u\ssp\|^{\,2}\ssp
   \big(
\sum_{\,i\ssp\in\sp\ssmb N\ssp\setminus\sp 2.}
i^{\,3.}\ssp\smb M^{\,i\yydot-\ssp 2}
+
\sum_{\,i\ssp=\sp\ssmb N}^{\ \infty}
i^{\,3.}\ssp 4^{\sp\,2\sp-\sp i\yydot}\sp
\big)$\vskip.7mm$\mhyppy{17.5}
\le|\ssp t\ssp|\,\|\ssp u\ssp\|^{\,2}\ssp
   \big(\sp
\smb N^{\,3.}\sp\smb M^{\,\ssmb N}\snn
+8\,
\sum_{\,i\ssp=\sp\ssmb N}^{\ \infty}
i^{\,2.}\ssp 2^{\sp\mminus i}\ssn\cdot
i\,2^{\sp\,1\sp-\ssp i\yydot}\sp
\big)$\vskip.7mm$\mhyppy{17.5}
\le|\ssp t\ssp|\,\|\ssp u\ssp\|^{\,2}\ssp(\sp
   \smb N^{\,3.}\sp\smb M^{\,\ssmb N}\snn
 + 8\sn\cdot\sn(\sp\ln 2\sp) ^{\sp\mminus 2}\ssn\cdot\sn
   2^{\sp\,2\ssp-\sp\ssmb N\yydot}(\spp 1+\smb N\ydot\sp))$\vskip.7mm$\mhyppy{17.5}
\le|\ssp t\ssp|\,\|\ssp u\ssp\|^{\,2}\ssp(\sp
   \smb N^{\,3.}\sp\smb M^{\,\ssmb N}\snn
 + 2^{\sp\,7\ssp-\sp\ssmb N\yydot}(\spp 1+\smb N\ydot\sp))\,$,\vskip.7mm

\noin whence we see that $\,\delta=\min\ssp\{\ssp 1\ssp,\sp
  \smb M^{\sp\mminus 2}\sp(\sp\smb N^{\,3.}\sp\smb M^{\,\ssmb N}\snn
+ 2^{\,7}\sp(\spp 1+\smb N\ydot\sp))^{\sp\mminus 1}\sp\eps\,\}\,$
will do.                                                       \end{example}


\affiliationone{Seppo I.\ Hiltunen\\
                 Helsinki University of Technology,
                 Institute of Mathematics, U311,
                 P.O.\ Box 1100,
                 FIN-02015 HUT\\
              FINLAND
   \email{shiltune\,@\,cc.hut.fi}}
\end{document}